\documentclass[a4paper,12pt,twoside]{article}
\usepackage[francais]{babel}
\usepackage{amssymb, amsmath, eucal}
\usepackage[all]{xy}
\usepackage{mathrsfs}
\begin{document}
\textwidth15.5cm
\textheight22.5cm
\voffset=-13mm
\newtheorem{The}{Theorem}[subsection]
\newtheorem{Lem}[The]{Lemma}
\newtheorem{Prop}[The]{Proposition}
\newtheorem{Cor}[The]{Corollary}
\newtheorem{Rem}[The]{Remark}
\newtheorem{Titre}[The]{\!\!\!\! }
\newtheorem{Conj}[The]{Conjecture}
\newcommand{\C}{\mathbb{C}}
\newcommand{\R}{\mathbb{R}}

\begin{center}

{\large\bf A Simple Proof of a Theorem by Uhlenbeck and Yau}

\end{center}

\vspace{3ex}

\begin{center}  {\bf Dan Popovici}    \end{center}

 {\small Universit\'e Joseph Fourier, Grenoble; Universit\'e de Paris-Sud, Orsay} 

\vspace{3ex}

\noindent {\bf Abstract.} {\small A subbundle of a Hermitian vector bundle $(E, \, h)$ can be metrically and differentiably defined by the orthogonal projection onto this subbundle. A weakly holomorphic subbundle of a Hermitian holomorphic bundle is, by definition, an orthogonal projection $\pi$ lying in the Sobolev space $L^2_1$ of $L^2$ sections with $L^2$ first order derivatives in the sense of distributions, which satisfies furthermore $(\mathrm{Id}-\pi)\circ D''\pi=0$. We give a new simple proof of the fact that a weakly holomorphic subbundle of $(E, \, h)$ defines a coherent subsheaf of ${\cal O}(E),$ that is a holomorphic subbundle of $E$ in the complement of an analytic set of codimension $\geq 2.$ This result was the crucial technical argument in Uhlenbeck's and Yau's proof of the Kobayashi-Hitchin correspondence on compact K\"ahler manifolds. We give here a much simpler proof based on current theory. The idea is to construct local meromorphic sections of $\mathrm{Im}\, \pi$ which locally span the fibers. We first make this construction on every one-dimensional submanifold of $X$ and subsequently extend it via a Hartogs-type theorem of Shiffman's.}

\vspace{3ex}

\subsection{Introduction}

 Let $(E, \,h)$ be a holomorphic vector bundle of rank $r$ equipped with a $C^{\infty}$ Hermitian metric over a compact K\"ahler manifold $X,$ and let ${\cal F}\subset {\cal O}(E)$ be a coherent analytic subsheaf of the locally free sheaf ${\cal O}(E)$ associated to $E.$ Since ${\cal F}$ is torsion-free (as a coherent subsheaf of a torsion-free sheaf), it is locally free outside an analytic subset of codimension $\geq 2$ (see, for instance, [Kob87], V.5). Thus ${\cal F}$ can be seen as a vector bundle with singularities. More precisely, there exists an analytic subset $S \subset X$, $\mathrm{codim}\, S \geq 2,$ and a holomorphic vector bundle $F$ on $X\setminus S,$ such that

\vspace{1ex}

\hspace{6ex}  ${\cal F}_{|X\setminus S} = {\cal O}(F).$ 

\vspace{1ex}

\noindent Equip the subbundle $F\hookrightarrow E_{|X\setminus S}$ with the Hermitian metric induced by $h$ and consider the orthogonal projection $\pi : E_{|X\setminus S} \longrightarrow F.$ Then $\pi$ can be seen as a  $C^{\infty}$ section on $X\setminus S$ of the holomorphic vector bundle $\mathrm{End}\, E$ satisfying: \\

$(0. 1)$  \hspace{2ex} $\pi = \pi^{\star} = \pi^2$, \hspace{3ex} $(\mathrm{Id}-\pi) \circ D''\pi = 0$ \\

\noindent on $X\setminus S,$ where $D''$ is the $(0, \, 1)$-component of the Chern connection on $\mathrm{End}\, E$ associated to the metric induced by $h$. The latter equality above says that the holomorphic structure of $F$ is the restriction of the holomorphic structure of $E_{|X\setminus S}.$ Let $Q$ be the quotient bundle of $E_{|X\setminus S}$ by $F,$ equipped with the metric induced by $h,$ and let $\det Q$ be the associated determinant line bundle equipped with the induced metric. Its curvature form $i\Theta(\det Q) = \mathrm{Tr}_Q(i\Theta(Q)) = \mathrm{Tr}_E(i\Theta(Q))$ is a $C^{\infty}$  $(1, \, 1)$-form on $X\setminus S$ given by the formula:

\vspace{1ex}

$i\Theta(\det Q) = \mathrm{Tr}_E(i\Theta_h(E)_{|Q}) + \mathrm{Tr}_E(i D'\pi \wedge D''\pi)$   \hspace{3ex}  (see [Gri69]).   

\vspace{1ex}

\noindent As $\mathrm{codim}\, S \geq 2,$ the $C^{\infty}$ $(1, \, 1)$-form  $\mathrm{Tr}_E(i D'\pi \wedge D''\pi)$ has locally finite mass in the neighbourhood of $S$. In other words, every $x\in S$ has a neighbourhood $U \subset X$ such that  \\

$\displaystyle\int_U \mathrm{Tr}_E(i D'\pi \wedge D''\pi) \wedge \omega^{n-1} < +\infty,$ \\

\noindent where $\omega$ is an arbitrary Hermitian metric on $X$. This is a consequence of a general current theory result stating that if $T$ is a closed positive current of bidegree $(p, \, p)$ (or equivalently of bidimension $(n-p, \, n-p)$) in the complement of an analytic subset $A$ of codimension $\geq p + 1$, then the mass of $T$ is locally finite in the neighbourhood of $A$ (see [Sib85], p. 178, {\it corollaire} 3. 2).

 In particular, the $(1, \, 1)$-form $\mathrm{Tr}_E(i D'\pi \wedge D''\pi)$ extended by $0$ across $S$ is $L^1$ on $X$. Since $|\mathrm{Tr}_E(i D'\pi \wedge D''\pi)|$ dominates  $|D'\pi|^2$ and  $|D''\pi|^2$, the norms being considered in their respective bundles, we see that $D'\pi$ and $D''\pi$ are $L^2$ $1$-forms on $X\setminus S$. Since every projection is $L^{\infty}$ and, thanks to the compacity of $X$, implicitely $L^2$, the projection $\pi$ belongs to the Sobolev space $L^2_1$ of $L^2$ sections of $\mathrm{End}\, E$ whose first order derivatives in the sense of distributions are still $L^2$. 

 This discussion can be summed up as follows. \\

\noindent {\bf Remark.} {\it Every coherent analytic subsheaf ${\cal F}$ of ${\cal O}(E)$ defines a section $\pi \in L^2_1(X, \, \mathrm{End}\, E)$ that is $C^{\infty}$ in the complement of an analytic set of codimension $\geq 2$ and satisfies relations $(0. 1)$.} \\

  The goal of the present paper is to prove, by relatively elementary methods, the reverse statement that was originally stated and proved in [UY 86, 89]. More precisely, we prove the following.

\begin{The}\label{The:principal} Let $(E, \, h)$ be a holomorphic vector bundle of rank $r$ equipped with a $C^{\infty}$ Hermitian metric over a compact complex K\"ahler manifold $X,$ and let $\pi \in L^2_1(X, \, \mbox{End}\, E)$ such that $\pi = \pi^{\star} = \pi^2$ and $(\mbox{Id}_E - \pi)\circ D''_{\mathrm{End} E}\pi = 0$ almost everywhere.

\vspace{2ex}

  Then there exist a coherent analytic subsheaf ${\cal F} \subset {\cal O}(E)$ and an analytic subset $S\subset X$ of codimension $\geq 2$ such that:

\vspace{2ex}

1) $\pi_{|X\setminus S}\in C^{\infty}(X\setminus S, \mbox{End}\,E)$

\vspace{2ex}

2) $\pi = \pi^{\star}= \pi^2$ and $(\mbox{Id}_E-\pi)\circ D''_{\mathrm{End}\, E}\pi = 0$ on $X\setminus S$

\vspace{2ex}

3) ${\cal F}_{|X\setminus S} = \pi_{|X\setminus S}(E_{|X\setminus S}) \hookrightarrow E_{|X\setminus S}$ is a holomorphic subbundle of $E_{|X\setminus S}$. 

\end{The}

\vspace{3ex}

 In what follows, $L^2_1(X, \, \mbox{End}\, E)$ stands for the Sobolev space of $L^2$ sections of the holomorphic bundle $\mbox{End}\, E$ whose first order derivatives in the sense of distributions are still $L^2.$ We equip $\mbox{End}\, E$ with the metric induced by $h$ and denote $D'_{\mathrm{End}\, E},$ $D''_{\mathrm{End}\, E}$ the $(1, \, 0)$ and respectively $(0, 1)$ components of the associated Chern connection. 

 A section $\pi \in L^2_1(X, \, \mbox{End}\, E)$ satisfying the hypotheses of theorem \ref{The:principal} is called {\it weakly holomorphic subbundle} of $E$.

  Theorem \ref{The:principal} provides the crucial technical argument in the proof given by Uhlenbeck and Yau ([UY 86, 89]) to the existence of a unique Hermitian-Einstein metric in every stable holomorphic vector bundle over a compact K\"ahler manifold. The comparatively easier converse, asserting that every Hermitian-Einstein vector bundle is semistable and splits into a direct sum of stable subbundles, had previously been proved by Kobayashi and L\"ubke ([Kob87, LT95]). Uhlenbeck and Yau were thus completing the proof of the Kobayashi-Hitchin correspondence over compact K\"ahler manifolds. The idea of their proof is the following. Having fixed the metric $h$ of $E,$ every $C^{\infty}$ metric $h_1$ on $E$ is of the form: 

\begin{center} $h_1(s, \, t) = h(f(s), \, t)$, \end{center}

\noindent for all sections $s$ and $t$ of $E$, where $f\in C^{\infty}(X, \, \mbox{End}\, E)$ is a positive definite and self-adjoint endomorphism (for $h$) of $E$. Then $h_1$ is a Hermitian-Einstein metric if and only if $f$ is a solution of a nonlinear partial differential equation. The authors solve a perturbed equation depending on a parameter $\varepsilon$ and find a solution $f_{\varepsilon}$. One of the following two situations occurs. Either $f_{\varepsilon}$ converges to an endomorphism $f_0$ when $\varepsilon$ tends to $0$, in which case they prove that $f_0$ actually defines a Hermitian-Einstein metric; or $f_{\varepsilon}$ does not converge, in which case they prove that the stability hypothesis on $E$ is violated by producing a destabilizing subsheaf of ${\cal O}(E)$. It is in the construction of such a destabilizing subsheaf that theorem \ref{The:principal} plays a key role.  

 Their proof is, however, extremely technical and not very enlightening. This is why we wish to give a natural proof of theorem \ref{The:principal} by constructing local meromorphic sections which span $\mathrm{Im}\, \pi$. The coherent sheaf ${\cal F}$ that we intend to construct will then be defined by its local sections.

\subsection{Preliminaries} Let $(X, \, \omega)$ be a compact K\"ahler manifold of dimension $n$ and let $(E, \, h)$ be a Hermitian holomorphic vector bundle of rank $r$ over $X.$ Let $\pi \in L^2_1(X, \, \mathrm{End}\, E)$ be a section such that $\pi=\pi^{\star}=\pi^2$ and $(\mathrm{Id}-\pi)\circ D''\pi=0$ almost everywhere, the derivative $D''\pi$ being computed in the sense of distributions. Since $\pi$ is a projection, we even have  

\vspace{1ex}

\hspace{6ex}  $\pi \in L^2_1(X,\, \mathrm{End}\, E) \cap L^{\infty}(X,\, \mathrm{End}\, E)$.

\vspace{1ex}

\noindent The subbundle $F=\mathrm{Im}\, \pi\subset E$ is defined almost everywhere as an $L^2$ bundle. This means that the fiber $F_x$ is defined as $\mathrm{Im}\, \pi_x$ for almost all points $x\in X$ and the transition matrices have an $L^2$ dependence on $x$. Likewise, the quotient bundle $Q=E/F$ is defined almost everywhere as an $L^2$ bundle. Let $\beta$ and $\beta^{\star}$ be the $(1, 0)$ $L^2$ current with values in $\mathrm{Hom}\, (F, \, Q)$, and respectively the $(0, 1)$ $L^2$ current with values in $\mathrm{Hom}\, (Q, \, F)$, uniquely determined by the following equalities:

\begin{center}

$D'\pi = \left( \begin{array}{cc}

                                   0 & 0 \\

                                   \beta & 0

                                        \end{array}  \right),$  \hspace{2ex} $D''\pi = \left( \begin{array}{cc}

                            0 & \beta^{\star}  \\

                            0 & 0

                 \end{array} \right), $  

\end{center}

\noindent where $D'\pi$ and $D''\pi$ are calculated in the sense of distributions. The current $\beta$ corresponds to the second fundamental form of the exact sequence $0\rightarrow F \rightarrow E \rightarrow Q \rightarrow 0$ in the case where $\pi$ is $C^{\infty}.$ We refer to [Gri69] for details concerning exact sequences of Hermitian vector bundles or, for a presentation using the same notation as in the present article, to chapter 5 of the book [Dem97].

\begin{Rem}\label{Rem:propdepibis} For all $\pi \in L^2_1(X, \, \mathrm{End}\, E)$ satisfying $\pi = \pi^{\star}=\pi^2$ almost everywhere, the following equalities in the sense of currents are equivalent: \\

$\begin{array}{lll}

(2. \, a) \,\,\,\,  (\mathrm{Id}-\pi)\circ D''\pi = 0 \, ;  & (2. \, b) &  D'\pi\circ (\mathrm{Id}-\pi)=0  \\

 (2. \, c) \,\,\,\,  \pi \circ D'\pi = 0 \, ;    & (2. \, d) &  D''\pi \circ \pi = 0.

\end{array}$  

\end{Rem}

\noindent {\it Proof.} The equivalence of $(a)$ and $(b)$ is obtained by taking adjoints, for $\pi=\pi^{\star}$. On the other hand, if we apply $D'$ to the equality $\pi=\pi^2$ we find $D'\pi = D'\pi \circ \pi + \pi \circ D'\pi$, which gives the equivalence of $(b)$ and $(c)$. Equality $(d)$ is inferred from $(c)$ by taking adjoints. All products are well-defined in the sense of currents since an $L^2$ form can be multiplied by an $L^{\infty}$ form. \hfill $\Box $  \\

 Proving theorem \ref{The:principal} amounts to proving that the $L^2$ bundle $F=\mathrm{Im}\, \pi$ is holomorphic outside an analytic subset of codimension $\geq 2.$ It is thus enough to construct local meromorphic sections of $F$ which span $F$ locally (for meromorphic sections are holomorphic in the complement of an analytic subset of codimension $\geq 2$). The idea is to construct local holomorphic sections of $F\otimes \det Q$ which span $F\otimes \det Q$ locally, in parallel with the construction of a local $\bar{\partial}$-closed section of $\det Q$ which spans $\det Q$ locally. A division of the local holomorphic sections of $F\otimes \det Q$ by the local section of $\det Q$ yields the meromorphic sections of $F$ that we wish to construct.

 In the construction of a local holomorphic section of $\det Q$ the $(1, 1)$-current $\mathrm{Tr}_E(i\beta\wedge \beta^{\star} + i\Theta(E)_{|Q})$, interpreted a posteriori as the curvature current of $\det Q,$ plays a key role.

\begin{Rem}\label{Rem:fermetureunevariable} The restriction of the current $\mathrm{Tr}_E(i\beta\wedge \beta^{\star} + (\mathrm{Id}-\pi) \circ i\Theta(E)_h \circ (\mathrm{Id}-\pi))$ to almost every complex line contained in a coordinate patch of $X$ defines a $d$-closed current.  

\end{Rem}

\noindent{\it Proof.} The argument is almost trivial. The existence of the restriction of an $L^1$ function to almost every line is a consequence of the Fubini theorem (the restriction being $L^1$ on this line). To see this, we start by considering a system of lines parallel to a given direction. Now, every current of maximal bidegree is closed. In particular, bidegree $(1, \, 1)$ currents are closed on complex submanifolds of dimension $1$.  \hfill $\Box$ \\

 In vue of theorem \ref{The:principal}, since the problem is local, we can work on an open set $U\subset X$ such that $E_{|U}\simeq U\times \C^r$. After possibly shrinking $U$, the curvature of $E$ can be made positive on $U$ by a change of metric. Indeed, let  

\vspace{1ex}

\hspace{6ex}  $h_1(z)=h(z) \cdot e^{-m|z|^2}$ 

\vspace{1ex}

\noindent be a new metric on $E_{|U}$, $m$ being a positive scalar and $z=(z_1, \dots , \, z_n)$ local coordinates on $U$. Since  \\

$i\Theta_{h_1}(E)=i\Theta_{h}(E) + m\, id'd''|z|^2\otimes \mathrm{Id}_E$,

\vspace{2ex}

\noindent and since the $(1, 1)$-form $id'd''|z|^2$ is $>0,$ we see that $i\Theta_{h_1}(E)\geq \varepsilon \omega \otimes \mathrm{Id}_{E}$ for a certain $\varepsilon > 0$, provided that $m$ is sufficiently large. On the other hand, the product $\beta \wedge \beta^{\star}$ defines an $L^1$ $(1, \, 1)$-current with values in $\mathrm{End}\, E$ and $\mathrm{Tr}_E(i\beta\wedge \beta^{\star}) \geq 0$ in the sense of currents (see [Gri69] for the $C^{\infty}$ case and the proof given there still works for currents). Thus the $(1, 1)$-current 

 $$\mathrm{Tr}_E(i\beta \wedge \beta^{\star} + (\mathrm{Id}-\pi)\circ i\Theta_{h_1}(E)\circ (\mathrm{Id}-\pi))$$

\noindent is positive on $U$. Moreover, this scalar change of metric preserves the property of $\pi$ being self-adjoint. We can then make the following convention without loss of generality. \\

\noindent {\bf Convention.} {\it We assume from now on that, locally, the curvature of $E$ is positive}. \\

\noindent Hence we get the following

\begin{Cor}\label{Cor:potentiel} The current $ \mathrm{Tr}_E(i\beta \wedge \beta^{\star} + (\mathrm{Id}-\pi)\circ i\Theta_{h_1}(E)\circ (\mathrm{Id}-\pi))$ of bidegree $(1, \, 1)$   admits a local subharmonic potential on almost every complex line contained in a coordinate patch of $X$. 

  This means that for every point $x \in X$ and for almost every complex line $L$ with respect to a system of local coordinates in the neighbourhood of $x$, there exists a subharmonic function $\varphi_L$ such that  \\ 

  $i\partial \bar{\partial}\varphi_L = \mathrm{Tr}_E(i\beta \wedge \beta^{\star} + (\mathrm{Id}-\pi)\circ i\Theta_{h_1}(E)\circ (\mathrm{Id}-\pi))$, \\

\noindent locally on $L$.

\end{Cor}

\noindent {\it Proof.} By the Poincar\'e lemma, every $d$-closed current is locally $d$-exact. It is therefore also $\partial \bar{\partial}$-exact by the $\partial \bar{\partial}$ lemma. After possibly shrinking the trivializing open set $U$, we may assume that there exists a function $\varphi_L$ as in the statement of the corollary. Since the above current is positive, the potential $\varphi_L$ is subharmonic.   \hfill  $\Box$

\vspace{2ex}

 The main difficulty in the proof of theorem \ref{The:principal} comes from the insufficient regularity of $\pi$. Certain wedge-products of currents are not well-defined for distributions cannot be multiplied. The following lemma states an elementary distribution theory result that will enable subsequent computations to make sense. For every real number $s,$ $L^2_s$ stands for the space of temperate distributions $u\in \mathscr{S}'(\R^n)$ such that the Fourier transform $\hat{u} \in L_{\mathrm{loc}}^1(\R^n)$ and $(1+|\xi|^2)^{\frac{s}{2}} \cdot \hat{u}(\xi) \in L^2(\R^n).$ It can be easily seen that if $u$ is a compactly supported element of $L^2_{s,\, \mathrm{loc}}$, there exists a decomposition $u=\sum\limits_j D^jv_j + v$, where $v_j\in L^2_{s+1,\, \mathrm{loc}}$ and $v\in L^2_{s+2, \,\mathrm{loc}}$ are compactly supported. In particular, for $s=-1$, every $L^2_{-1}$ function can be locally written as a sum of partial derivatives of order 1, in the sense of distributions, of $L^2$ functions, modulo an $L^2_1$ function. In an analogous way we define the space $(L^1)'$ of temperate distributions arising locally as a sum of partial derivatives of order 1, in the sense of distributions, of $L^1$ functions. We thus have the inclusions : 

\vspace{1ex}

\hspace{6ex}  $L^1 \hookrightarrow (L^1)' \hookrightarrow {\cal D}'_1$, 

\vspace{1ex}

\noindent where ${\cal D}'_1$ stands for the space of order 1 distributions. The topology of $(L^1)'$ is defined to be the restriction of the topology of ${\cal D}'_1$.

 \begin{Lem}\label{Lem:distr}

 The map

\vspace{1ex}

$(f, \, g)\mapsto u_{fg}$ \hspace{2ex} acting from $L^2_{1,\, \mathrm{loc}}\times L^2_{-1, \, \mathrm{loc}}$ to $(L^1)'$ 

\vspace{1ex}

\noindent is well-defined, bilinear and continuous, where $u_{fg}$ stands for the distribution defined as:

\vspace{1ex}

$\displaystyle <u_{fg}, \, \varphi>=-\sum\limits_j\int g_jD^j(\theta_j f\varphi) + \int f  \psi h \varphi,  $ 

\vspace{1ex}

\noindent for every test function $\varphi$ and every local decomposition $g=\sum\limits_j \theta_j\cdot D^j g_j + \psi h $, with $g_j \in L^2, \, h\in L^2_1$ and $\theta_j, \, \psi$ test functions.

\end{Lem}

 The standard proof of this lemma can well be left to the reader.

\subsection{Proof of Theorem \ref{The:principal}}

 The overall idea of proof is to achieve enough regularity on expressions containing $\pi$ and derivatives of $\pi$ enabling us to retrieve the classical $C^{\infty}$ situation outside an analytic set of codimension $\geq 2.$ Frequent side-glances at the $C^{\infty}$ situation will show us the way. We shall proceed in several steps. \\

\noindent $\bullet$ {\bf First step} : {\it reduction to the case of zero curvature}

 In order to massively simplify subsequent computations, we start off by showing that we can locally reduce the problem to the case where the curvature of $E$ vanishes identically. The following elementary lemma will be of use.

\begin{Lem}\label{Lem:metriques} Let $E$ be a complex vector space of dimension $r$ and $F$ a vector subspace of dimension $p$. Consider two Hermitian metrics $h$ and $h_0$ on $E$ and let $\pi$, $\pi_0$ be the orthogonal projections, for $h$ and respectively $h_0,$ of $E$ on $F$. 

  If $E= F \oplus F_h^{\perp}$ (respectively $E= F \oplus F_{h_0}^{\perp}$) is the orthogonal decomposition of $E$ for $h$ (respectively  $h_0$), then there exists an automorphism $v: E \longrightarrow E$ such that $v(F)=F$, $v(F_h^{\perp}) = F_{h_0}^{\perp}$, and $h(s, \, t) = h_0(vs, \, vt)$, for all $s, t \in E$.

  Moreover, for every such $v$ the projections $\pi$ and $\pi_0$ are related by the formula $\pi_0 = v\pi v^{-1}$.

\end{Lem}

 The elementary proof of this lemma is left to the reader.

\begin{Cor}\label{Cor:metriques} Let $(E, \, h)$ be a holomorphic vector bundle of rank $r$ equipped with a $C^{\infty}$ Hermitian metric over a complex manifold $X,$ and let $\pi \in L^2_1(X, \, \mathrm{End}\, E)$ be such that $\pi = \pi^{\star} = \pi^2$ and $(\mathrm{Id}_E - \pi) \circ D''_{\mathrm{End}\, E} \pi = 0$ almost everywhere. Set $F = \mathrm{Im}\, \pi$. Let $U$ be a trivializing open set for the bundle $E$ and let $h_0$ be the trivial flat metric on $E_{|U}\simeq U \times \C^r$. Let $\pi_0 \in L^2_1(U, \, \mathrm{End}\, E)$ be the orthogonal projection of $E_{|U}$ onto $F_{|U}$ with respect to the metric $h_0$.

  Then there exists $v\in C^{\infty}(U, \, \mathrm{End}\, E)$ such that $(\mathrm{Id}-\pi)\circ v \circ \pi = 0$ (or equivalently $(\mathrm{Id}-\pi_0)\circ v \circ \pi_0 = 0$), $\pi_0 \circ v \circ (\mathrm{Id}-\pi)=0$ almost everywhere on $U,$ and $h(s, \, t) = h_0(vs, \, vt)$ for all $s, \, t \in E_{|U}$. Furthermore, $\pi_0 = v\pi v^{-1}$ almost everywhere on $U$.

\end{Cor}

\begin{Lem}\label{Lem:pizero} Under the hypotheses of corollary \ref{Cor:metriques} the projection $\pi_0$ satisfies moreover: $(\mathrm{Id}-\pi_0) \circ D''\pi_0 = 0$ almost everywhere on $U$.

\end{Lem}

 This result corresponds a posteriori to the holomorphic structure of $F,$ viewed as a holomorphic subbundle of $E$ in the complement of an analytic set, being independent of the choice of metric. \\

\noindent {\it Proof.} As $\pi_0 = v\pi v^{-1}$, we infer: \\

$\begin{array}{lll}

(\mathrm{Id}-\pi_0)\circ D''\pi_0 & = & (\mathrm{Id}-\pi_0)\circ D''v \circ \pi \circ v^{-1} \\

  & + & (\mathrm{Id}-v\pi v^{-1}) \circ v \circ D''\pi \circ v^{-1} \\

  & + & (\mathrm{Id}-v\pi v^{-1}) \circ v \circ \pi \circ D''(v^{-1}).

\end{array}$ \\

\noindent The above expressions are well-defined in the sense of distributions since $v$ is $C^{\infty}$. The second term in the above sum is equal to: \\

$v\circ D''\pi \circ v^{-1} - v\circ \pi \circ D''\pi \circ v^{-1} = 0$, \\

\noindent for $D''\pi = \pi \circ D''\pi$ and the subtraction of two $L^2$ expressions is well-defined. The third term in the above sum is equal to: \\

$v \circ \pi \circ D''(v^{-1}) - v \circ \pi \circ v^{-1} \circ v \circ \pi   \circ D''(v^{-1}) = 0$, \\

\noindent for $\pi \circ v^{-1} \circ v \circ \pi = \pi^2 = \pi$ and two $L^2$ expressions can be subtracted. The sum is thus reduced to its first term. Hence: \\

 $(a)$  \hspace{2ex} $ (\mathrm{Id}-\pi_0)\circ D''\pi_0 =    (\mathrm{Id}-\pi_0)\circ D''v \circ \pi \circ v^{-1}$. \\

\noindent Apply the operator $D''$ to the equality $(\mathrm{Id}-\pi) \circ v \circ \pi = 0$ and get: \\

$(b)$  \hspace{2ex} $-D''\pi \circ v \circ \pi + (\mathrm{Id}-\pi) \circ D''v \circ \pi + (\mathrm{Id}-\pi) \circ v \circ D''\pi = 0$. \\

\noindent Since $D''\pi = \pi \circ D''\pi$, we see that: \\

  $(\mathrm{Id}-\pi) \circ v \circ D''\pi = \bigg ((\mathrm{Id}-\pi) \circ v \circ \pi \bigg )\circ D''\pi = 0$, \\

\noindent for the expression in brackets vanishes. Equality $(b)$ becomes: \\

$(c)$  \hspace{2ex} $D''\pi \circ v \circ \pi = (\mathrm{Id}-\pi) \circ D''v \circ \pi$. \\

\noindent On the other hand, for all $\xi\in E$ there exists $\eta \in E$ such that $v(\pi \xi) = \pi \eta$, for $v$ preserves $\mathrm{Im}\, \pi$. Since $D''\pi \circ \pi = 0$ thanks to relation $(2. \, d)$, we get: \\

$(D''\pi \circ v \circ \pi)(\xi) = D''\pi (v(\pi \xi)) = (D''\pi \circ \pi) (\eta) = 0$, \\

\noindent for all $\xi \in E$. Consequently, $D''\pi \circ v \circ \pi = 0$ and relation $(c)$ implies: $(\mathrm{Id}-\pi) \circ D''v \circ \pi = 0$. This is equivalent to: $D''v \circ \pi = \pi \circ D''v \circ \pi$. If we apply the operator $\mathrm{Id}-\pi_0$ to this last equality we get: \\

$(\mathrm{Id}-\pi_0) \circ D''v \circ \pi = (\mathrm{Id}-\pi_0) \circ \pi \circ D''v \circ \pi = 0$, \\

\noindent for $(\mathrm{Id}-\pi_0) \circ \pi = 0$ ( $\mathrm{Im}\, \pi = \mathrm{Im}\, \pi_0$ and $(\mathrm{Id}-\pi_0) \circ \pi_0 = 0$). Equality $(a)$ finally gives: $(\mathrm{Id}-\pi_0) \circ D''\pi_0 = 0$, and this is what we wanted to prove.   \hfill $\Box$ \\

 Lemma \ref{Lem:pizero} enables us to reduce locally to the case of a flat vector bundle. Indeed, since the problem is local we may assume from now on, possibly after replacing locally the original metric $h$ of $E$ with the trivial flat metric $h_0$, that $i\Theta(E)_h = 0$ on the trivializing open set $U$.

 \vspace{2ex}

\noindent $\bullet$ {\bf Second step} : {\it reinterpretation of $\mathrm{Im}\, \pi$}

  In order to prove theorem \ref{The:principal} we will show that the $L^2$ bundle $F=\mathrm{Im}\, \pi$ is locally generated by its local meromorphic sections. As before, we will draw on the $C^{\infty}$ situation considered in the very simple lemma below.

\begin{Lem}\label{Lem:morphismesurj} Let $(E, \, h)$ be a Hermitian holomorphic vector bundle of rank $r$ and let $\pi\in C^{\infty}(X, \, \mathrm{End}\, E)$ be such that $\pi = \pi^{\star}=\pi^2$ and $(\mathrm{Id}-\pi)\circ D''\pi = 0$. We assume that the curvature form of $(E, \, h)$ vanishes identically. Let $p$ be the rank of $\pi$ and $q=r-p$. Consider the holomorphic subbundle $F=\mathrm{Im}\, \pi$ of $E$ and the exact sequence of holomorphic vector bundles: \\

$0 \longrightarrow F \stackrel{j}{\longrightarrow} E \stackrel{g}{\longrightarrow} Q \longrightarrow 0$, \\

\noindent where $j$ is the inclusion and $g= \mathrm{Id}-\pi$ is the projection onto the quotient bundle $Q$. Then there exists a holomorphic bundle morphism: \\

$\Lambda^{q+1}E \otimes \Lambda^q Q^{\star}   \stackrel{\sigma}{\longrightarrow} E, $ \\

\noindent whose image is $F$. More precisely, if $e_1, \dots , \, e_r$ is a local orthonormal holomorphic frame of $E$ and $K=(k_1 < \dots < k_q)$ is a multiindex, consider the local holomorphic section of $\det Q = \Lambda^qQ$ definied as: \\

$v_K = (\mathrm{Id}-\pi)e_{k_1} \wedge \dots \wedge (\mathrm{Id}-\pi)e_{k_q} = \sum\limits_{|J|=q}D_{JK} \cdot e_J$ \\

\noindent where $D_{JK}$ is the minor corresponding to the lines $J=(j_1 < \dots < j_q)$ and the columns $K=(k_1 < \dots < k_q)$ in the matrix representing $\mathrm{Id}-\pi$ in the frame under consideration, and $e_J:=e_{j_1} \wedge \dots \wedge e_{j_q}$ for all $J=(j_1 < \dots < j_q)$. Associate to $v_K$ the local holomorphic section of $\Lambda^qQ^{\star}$ definied as: \\

$\displaystyle v_K^{-1} = \frac{\sum\limits_{|J|=q}\bar{D}_{JK} \cdot e_J^{\star}}{\sum\limits_{|J|=q}|D_{JK}|^2}.$ \\

  Then for all multiindices $I=(i_1 <  \dots <  i_{q+1})$ and $K=(k_1 < \dots < k_q)$, the morphism $\sigma$ is locally defined by the relation: \\

$(3. 1)$  \hspace{2ex}  $\displaystyle \sigma(e_I \otimes v_K^{-1})= \overset{q+1}{\underset{l=1}\sum} (-1)^l \cdot \frac{\sum\limits_{|J|=q}\bar{D}_{JK} \cdot e_J^{\star}(e_{I\setminus \{i_l\}})}{\sum\limits_{|J|=q}|D_{JK}|^2} \cdot e_{i_l}.$

\noindent  In particular, there exists a holomorphic bundle morphism: \\

$\Lambda^{q+1}E \stackrel{u}{\longrightarrow} E \otimes \det Q$ \\

\noindent whose image is $F\otimes \det Q$, which is induced by $\sigma$ after tensorizing to the right by $\det Q = \Lambda^qQ$. Morphisms $\sigma$ and $u$ are locally related by: \\

$\displaystyle\sigma(e_I \otimes v_K^{-1}) = \frac{u(e_I)}{v_K}$, \\

\noindent the division being performed in the line bundle $\det Q$.

\end{Lem}

\vspace{2ex}

 This lemma shows that the vector bundle $F\otimes \det Q$ can be realized as the image of a holomorphic projection from $\Lambda^{q+1}E$. It will prove useful later on since the projection of $E$ on $F$ is not holomorphic in general. \\

\noindent {\it Proof.} The quotient bundle $Q$ can be seen as a $C^{\infty}$ subbundle of $E$ via the $C^{\infty}$ inclusion $Q \stackrel{\mathrm{Id}-\pi}{\hookrightarrow}E.$ This defines the $C^{\infty}$ inclusion $\det Q = \Lambda^qQ \hookrightarrow \Lambda^qE$ and the orthogonal decomposition $\Lambda^qE = \Lambda^qQ \oplus (\Lambda^qQ)^{\perp}$. The element $v_K^{-1}$ of $\Lambda^q E^{\star}$ satisfies the identities: \\

$v_K^{-1}(v_K) = 1$  \hspace{2ex} and \hspace{2ex} $v_K^{-1}(\xi) = 0$, \, for all $\xi \in (\Lambda^qQ)^{\perp}$. \\

\noindent This accounts for the notation $v_K^{-1}$ and shows that $v_K^{-1}\in \Lambda^q Q^{\star}=(\det Q)^{-1}.$ We have to prove that $\mathrm{Im}\, \sigma = F$. The inclusion $F \subset \mathrm{Im}\, \sigma$ is obvious. Indeed, for all $s\in F$, $\sigma(s\wedge v_K \otimes v_K^{-1}) = v_K^{-1}(v_K) \cdot s = s$. Let us now prove the inclusion $\mathrm{Im}\, \sigma \subset F$. We have: $v_K^{-1} \in \Lambda^qQ^{\star} \hookrightarrow \Lambda^qE^{\star}$, and the inclusion is holomorphic. Viewed as an element of $\Lambda^qE^{\star}$, $v_K^{-1}$ is defined by: \\

$v_K^{-1}(e_{j_1} \wedge \dots \wedge e_{j_q}) = v_K^{-1}((\mathrm{Id}-\pi)e_{j_1} \wedge \dots \wedge (\mathrm{Id}-\pi)e_{j_q}),$  \\

\noindent for all $1 \leq j_1 < \dots < j_q \leq r$. Consequently, for all multiindex $I = (i_1 < \dots < i_{q+1})$ we have : \\

$\sigma(e_I \otimes v_K^{-1})  =  \overset{q+1}{\underset{l=1}\sum}(-1)^l \, v_K^{-1}(e_{i_1} \wedge \dots \wedge \hat{e}_{i_l} \wedge \dots \wedge e_{i_{q+1}}) \, e_{i_l}$ \\

$   =  \overset{q+1}{\underset{l=1}\sum}(-1)^l \, v_K^{-1}\bigg((\mathrm{Id}-\pi)e_{i_1} \wedge \dots \wedge \widehat{(\mathrm{Id}-\pi)e}_{i_l} \wedge \dots \wedge (\mathrm{Id}-\pi)e_{i_{q+1}}\bigg) \, e_{i_l}$ \\

$   =  \overset{q+1}{\underset{l=1}\sum}(-1)^l \, \sigma\bigg((\mathrm{Id}-\pi)e_{i_1} \wedge \dots \wedge e_{i_l} \wedge \dots \wedge (\mathrm{Id}-\pi)e_{i_{q+1}}\otimes v_K^{-1}\bigg)$ \\

$   =  \overset{q+1}{\underset{l=1}\sum}(-1)^l \, \sigma\bigg((\mathrm{Id}-\pi)e_{i_1} \wedge \dots \wedge \pi e_{i_l} \wedge \dots \wedge (\mathrm{Id}-\pi)e_{i_{q+1}}\otimes v_K^{-1}\bigg)$ \\

$  =  \overset{q+1}{\underset{l=1}\sum}(-1)^l \, v_K^{-1}\bigg((\mathrm{Id}-\pi)e_{i_1} \wedge \dots \wedge \widehat{(\mathrm{Id}-\pi)e}_{i_l} \wedge \dots \wedge \mathrm{Id}-\pi)e_{i_{q+1}}\bigg) \, \pi e_{i_l}.$

\noindent This proves that $\sigma(e_I \otimes v_K^{-1}) \in F$, for all multiindices $I$ and $K$, since $\pi e_{i_l}  \in F$ for all $i_l$. Therefore, $\mathrm{Im}\, \sigma \subset F$.   \hfill $\Box$ \\

\vspace{3ex}

 Lets us now turn back to the situation in theorem \ref{The:principal}. The section $\pi$ under consideration is $L^2_1$. Fix a local holomorphic frame $e_1, \dots , \, e_r$ of $E$ on an open set $U$ and let $p$ be, as above, the rank almost everywhere of $F=\mathrm {Im}\, \pi,$ and $q=r-p$. For a fixed point $x_0 \in U$ we may assume that $e_1(x_0), \dots , \, e_q(x_0)$ is a basis of $Q_{x_0}$ and that $e_{q+1}(x_0), \dots , \, e_r(x_0)$ is a basis of $F_{x_0}$. It is then obvious that $(\mathrm{Id}-\pi)e_j(x_0) = e_j(x_0)$ for $j\in \{1, \dots ,\, q\},$ and that $(\mathrm{Id}-\pi)e_j(x_0) = 0$ for $j\in \{q+1, \dots , \, r\}.$ For every matrix $a=(a_{k,  j})_{\stackrel{1\leq k \leq q}{1\leq j \leq r}}$ with $a_{k,  j}\in \C$ such that $(a_{k,  j})_{\stackrel{1\leq k \leq q}{1\leq j \leq q}} = \mathrm{Id}_{\C^q}$, let us define the following local holomorphic sections of $E$ on $U$: \\

$s_k=\overset{r}{\underset{j=1}\sum}a_{k,  j}e_j$, \, for $k=1, \dots , \, q,$ \\

\noindent and the local section of $\Lambda^qE$ on $U$: \\

 $\tau_a = (\mathrm{Id}-\pi)s_1 \wedge \dots \wedge (\mathrm{Id}-\pi)s_q \in L^2_1 \cap L^{\infty}$. \\

\noindent The section $\tau_a$ is a linear combination of the sections $v_K$ of $\det Q$ considered in lemma \ref{Lem:morphismesurj}. A posteriori $\tau_a$ will be a local holomorphic section of $\det Q$. Furthermore, $\tau_a(x_0) = e_1(x_0) \wedge \dots \wedge e_q(x_0)$, and therefore $|\tau_a(x_0)| \neq 0$.   Imitating formula $(3. 1)$ of lemma \ref{Lem:morphismesurj}, we thus get

\begin{Cor} For a Hermitian vector bundle $(E, \, h)$ and a section $\pi\in L^2_1(X, \, \mathrm{End}\, E)$ satisfying the hypotheses of theorem \ref{The:principal}, the local bundle morphism $\Lambda^{q+1}E_{|U} \stackrel{v}{\longrightarrow} E_{|U}$ defined by the formula: \\

$(3. 2)$  \hspace{2ex} $\displaystyle e_I:= e_{i_1}\wedge \dots \wedge e_{i_{q+1}} \stackrel{v}{\longmapsto}  \sigma(e_I \otimes \tau_a^{-1}) =  \frac{u(e_I)}{\tau_a},$ \\

\noindent for all $I=(1 \leq i_1 <  \dots  \, < i_{q+1} \leq r),$ satisfies $\mathrm{Im}\, \pi_{|U} = \mathrm{Im}\, v.$ Here $\sigma$ and $u$ are defined by the same formulae as in lemma \ref{Lem:morphismesurj}.

\end{Cor}

\noindent $\bullet$ {\bf Third step} : {\it a Lelong-Poincar\'e-type lemma}

 This is the key step of the proof. Since $\mathrm{Im}\, \pi =  \mathrm{Im}\, v$ locally, it is enough to prove that, for every multiindex $I$ such that $|I|=q+1$ we have $D''(v( e_I)) = 0$ in the sense of currents, in order to conclude that $\mathrm{Im}\, \pi$ defines a holomorphic bundle outside an analytic subset of codimension $\geq 2$. Although formula $D'' (v( e_I)) = 0$ is formally true, this makes a priori no sense for $\frac{1}{\tau_a}$ does not necessarily define a distribution. (The coefficients of $\tau_a$ are $L^2_1$ functions and their inverses are only measurable functions.) What is therefore at stake in this approach is to achieve enough regularity enabling us to apply the operator $D''$ in the sense of distributions.

 To begin with, let us notice that a posteriori the Lelong-Poincar\'e formula applied to the holomorphic section $\tau_a$ of the (a posteriori) holomorphic line bundle $\det Q$ gives: \\

$\frac{i}{2\pi}\partial\bar{\partial} \log |\tau_a| = [Z_a] - \frac{i}{2\pi} \Theta(\mbox{d\'et}\, Q) = [Z_a] - \frac{1}{2\pi} \mathrm{Tr}_E(i\beta\wedge \beta^{\star}),$ \\

\noindent where $[Z_a]$ stands for the current of integration along the zero divisor $Z_a$ of $\tau_a$ and $|\tau_a|$ designates the quotient norm of $\tau_a$ in $\det Q$ (equal, as a matter of fact, to the norm of $\tau_a$ in $\Lambda^qE$). Since the curvature of $E$ is assumed to vanish identically, $i\Theta(\det Q) = \mathrm{Tr}_E(i\beta\wedge \beta^{\star})$. Since $[Z_a]$ is a $(1, \, 1)$ positive current, we get: \\

\hspace{3ex} $i\partial\bar{\partial} \log |\tau_a|^2 \geq  - \mathrm{Tr}_E(i\beta\wedge \beta^{\star}).$ \\

\noindent This a posteriori situation can be retrieved by direct computation of $i\partial\bar{\partial} \log |\tau_a|^2$ (the norm being taken in $\Lambda^q E$). This is the goal of the present step of the proof. Consider the following operators: \\

\hspace{6ex} $D'_Q:= (\mathrm{Id} - \pi) \circ D'_E$ ,  \hspace{3ex}  $D''_Q:= (\mathrm{Id} - \pi) \circ D''_E$, \\

\noindent representing the projections of $D'_E$ et respectively $D''_E$. A posteriori, $D'_Q$ and $D''_Q$ will be the $(1, 0)$ and respectively $(0, 1)$ components of the Chern connection associated to the quotient metric of the vector bundle $Q$ that we will construct. In order to avoid complications with denominators we will compute $i\partial\bar{\partial}\log (|\tau_a|^2 + \delta^2)$ for real numbers $\delta > 0$ that we shall subsequently make converge to $0$. The following lemma yields the key argument to the proof of theorem \ref{The:principal}.

\begin{Lem}\label{Lem:LelPoinc} If $(E, \, h)$ is a holomorphic vector bundle of rank $r$ equipped with a $C^{\infty}$ Hermitian metric of zero curvature form, then for all $\delta > 0$ we have \\

$\begin{array}{lll}

\displaystyle i\partial\bar{\partial} \log(|\tau_a|^2 + \delta^2) & = &  i\frac{ \{D'_{\mbox{\tiny det\!\! Q}}\tau_a, \, D'_{\mbox{\tiny det Q}} \tau_a \} }{|\tau_a|^2 + \delta^2} - i \frac{ \{D'_{\mbox{\tiny det Q}}\tau_a, \, \tau_a\} \wedge \{\tau_a, \, D'_{\mbox{\tiny det Q}}\tau_a \}}{(|\tau_a|^2 + \delta^2)^2}  \\

\vspace{2ex} 

   & - & i\frac{\{D''_{\mbox{\tiny det Q}}D'_{\mbox{\tiny det\!\! Q}}\tau_a, \, \tau_a \} } {|\tau_a|^2 + \delta^2} \\

\vspace{2ex} 

 & \geq & -\frac{|\tau_a|^2}{|\tau_a|^2 + \delta^2} \cdot \mathrm{Tr}_E (i\beta \wedge \beta^{\star}),

\end{array}$ \\

\noindent where 

$D'_{\mbox{\tiny det\! Q}}\tau_a := \overset{q}{\underset{k=1}\sum}(\mathrm{Id}-\pi)s_1 \wedge \dots \wedge D'_Q\bigg((\mathrm{Id}-\pi)s_k\bigg) \wedge \dots \wedge (\mathrm{Id}-\pi)s_q,$  \\

$D''_{\mbox{\tiny det\! Q}}\tau_a := \overset{q}{\underset{k=1}\sum}(\mathrm{Id}-\pi)s_1 \wedge \dots \wedge D''_Q\bigg((\mathrm{Id}-\pi)s_k\bigg) \wedge \dots \wedge (\mathrm{Id}-\pi)s_q.$

\end{Lem}

 A word of explanation is in order here to justify the well-definedness of the above expressions. Recall that $\pi \in L^{\infty}(X, \, \mathrm{End}\, E) \cap L^2_1(X, \, \mathrm{End}\, E)$. The currents $\{D'_{\mbox{\tiny det\!\! Q}}\tau_a, \, D'_{\mbox{\tiny d\'et\!\! Q}} \tau_a \}$ and $\{D'_{\mbox{\tiny det\!\! Q}}\tau_a, \, \tau_a\} \wedge \{\tau_a, \, D'_{\mbox{\tiny det\!\! Q}}\tau_a \}$ are well-defined $(1, 1)$-currents with $L^1$ coefficients because $D'_{\mbox{\tiny det\!\! Q}} \tau_a$ is a current with $L^2$ coefficients (arising as products of an $L^2$ function by $q-1$ $L^{\infty}$ functions). On the other hand, $\frac{\tau_a}{|\tau_a|^2 + \delta^2}$ is easily seen to be an $L^2_1$ section of $\Lambda^qE$. Then lemma \ref{Lem:distr} implies that the $(1, 1)$-current : \\

$i\frac{\{D''_{\mbox{\tiny det\!\! Q}}D'_{\mbox{\tiny det\!\! Q}}\tau_a, \, \tau_a \} } {|\tau_a|^2 + \delta^2} = i\{D''_{\mbox{\tiny det\!\! Q}}D'_{\mbox{\tiny det\!\! Q}}\tau_a, \, \frac{\tau_a}{|\tau_a|^2 + \delta^2} \}$ \\

\noindent is well-defined and has $(L^1)'$ distributions as coefficients obtained as products of an $L^2_{-1}$ distribution by an $L^2_1$ function. \\

\noindent {\it Proof of lemma \ref{Lem:LelPoinc}.} Since the curvature form of $E$ is assumed to be zero (see the first step), the local frame $e_1, \dots , \, e_r$ can be chosen to be parallel for the Chern connection of $(E, \, h)$. This means that $D'_Ee_j=0$ and $D''_Ee_j=0$, for all $j=1, \dots , \, r$. The section $\tau_a$ of $\Lambda^qE$ can then be locally written with respect to this frame as: \\

$\tau_a = \sum\limits_{j_1, \dots , \, j_q} a_{1 j_1} \dots a_{qj_q} (\mathrm{Id}-\pi)e_{j_1} \wedge \dots \wedge (\mathrm{Id}-\pi)e_{j_q}.$  \\

\noindent The operators $D'_{\Lambda^qE}$ and $D'_{\mbox{\tiny det\!\! Q}}$, as well as the operators $D''_{\Lambda^qE}$ and $D''_{\mbox{\tiny det\!\! Q}}$, when applied to $\tau_a$ yield: \\

$D'_{\Lambda^qE}\tau_a = - \sum\limits_{j_1, \dots , \, j_q} a_{1 j_1} \dots a_{qj_q} \cdot \sum\limits_k (\mathrm{Id}-\pi)e_{j_1} \wedge \dots  \wedge D'\pi(e_{j_k}) \wedge    \dots \wedge (\mathrm{Id}-\pi)e_{j_q}; $ \\

$D'_{\mbox{\tiny det\!\! Q}}\tau_a =  - \sum\limits_{j_1, \dots , \, j_q} a_{1 j_1} \dots a_{qj_q} \cdot \sum\limits_k (\mathrm{Id}-\pi)e_{j_1} \wedge \dots  \wedge (\mathrm{Id}-\pi)\circ D'\pi(e_{j_k}) \wedge    \dots \wedge (\mathrm{Id}-\pi)e_{j_q}; $ \\

$D''_{\Lambda^qE}\tau_a =  - \sum\limits_{j_1, \dots , \, j_q} a_{1 j_1} \dots a_{qj_q} \cdot \sum\limits_k (\mathrm{Id}-\pi)e_{j_1} \wedge \dots  \wedge D''\pi(e_{j_k}) \wedge    \dots \wedge (\mathrm{Id}-\pi)e_{j_q}; $  \\

$D''_{\mbox{\tiny det\!\! Q}}\tau_a =  - \sum\limits_{j_1, \dots , \, j_q} a_{1 j_1} \dots a_{qj_q} \cdot \sum\limits_k (\mathrm{Id}-\pi)e_{j_1} \wedge \dots  \wedge (\mathrm{Id}-\pi)\circ D''\pi(e_{j_k}) \wedge    \dots \wedge (\mathrm{Id}-\pi)e_{j_q}. $ \\

\noindent The currents $D'_{\Lambda^qE}\tau_a,$ $D''_{\Lambda^qE}\tau_a,$ $D'_{\mbox{\tiny det\!\! Q}}\tau_a$ and $D''_{\mbox{\tiny det\!\! Q}}\tau_a$ are well-defined currents with $L^2$ coefficients. Since $(\mathrm{Id}-\pi)\circ D'\pi = D'\pi,$ as implied by relation $(2. \, c)$, we see that 

\vspace{1ex}

$D'_{\mbox{\tiny det\!\! Q}}\tau_a = D'_{\Lambda^qE}\tau_a.$

\vspace{1ex}

\noindent Let us denote from now on this common value by $D'\tau_a$. Relation $(2. \, a)$ shows that $(\mathrm{Id}-\pi)\circ D''\pi = 0$, which entails that $D''_{\mbox{\tiny det\!\! Q}}\tau_a = 0$.

 Let us start off by proving the inequality featured by lemma \ref{Lem:LelPoinc}. The $(1,\, 1)$-current $i\{D'\tau_a, \, D'\tau_a \}$ is positive and the Cauchy-Schwarz inequality shows that: \\

 $i\{D'\tau_a , \, \tau_a \} \wedge \{\tau_a, \, D'\tau_a \} \leq |\tau_a|^2 \cdot i\{D'\tau_a, \, D'\tau_a\}$, \\

\noindent in the sense of currents. Consequently, $i\frac{ \{D'\tau_a, \, D' \tau_a \}}{|\tau_a|^2 + \delta^2} - i \frac{ \{D'\tau_a, \, \tau_a\} \wedge \{\tau_a, \, D'\tau_a \}}{(|\tau_a|^2 + \delta^2)^2} \geq 0$, in the sense of currents. We will now show the following equality: \\

$({\star})$  \hspace{2ex} $ \{D''_{\mbox{\tiny det\!\! Q}}D'_{\mbox{\tiny det\!\! Q}}\tau_a, \, \tau_a \} = |\tau_a|^2  \cdot \mathrm{Tr}_E (\beta \wedge \beta^{\star}),$ \\

\noindent and this will prove the inequality stated in the lemma. Since: 

\vspace{2ex}

$D''_{\mbox{\tiny det\!\! Q}} ((\mathrm{Id}-\pi)e_j)= - (\mathrm{Id}-\pi)\circ D''\pi (e_j)=0$, \hspace{2ex} for all $j$, we get: \\

$ D''_{\mbox{\tiny det\!\! Q}} D'_{\mbox{\tiny det\!\! Q}}\tau_a = - \sum\limits_{j_1, \dots , \, j_q} a_{1 j_1} \dots a_{qj_q} $ \\

\hspace{6ex} $\cdot \sum\limits_k (\mathrm{Id}-\pi)e_{j_1} \wedge \dots  \wedge (\mathrm{Id}-\pi)\circ D''D'\pi(e_{j_k}) \wedge    \dots \wedge (\mathrm{Id}-\pi)e_{j_q}.$ \\

\noindent On the other hand, upon applying the operator $D'$ to the identity $(\mathrm{Id}-\pi)\circ D''\pi = 0$ we get: 

\vspace{2ex}

$-D'\pi \wedge D''\pi + (\mathrm{Id}-\pi)\circ D'D''\pi = 0,$ 

\vspace{2ex}

\noindent in the sense of currents. The current $D'\pi \wedge D''\pi$ has $L^1$ coefficients since it is the product of two currents with $L^2$ coefficients. The current $D'D''\pi$ has $L^2_{-1}$ coefficients and lemma \ref{Lem:distr} allows its multiplication by the current $\mathrm{Id}-\pi$ with $L^2_1$ coefficients. Furthermore, since the curvature of $E$ is assumed to be zero, the curvature of $\mathrm{End}\, E,$ equipped with the metric induced from the metric of $E,$ is zero as well. Thus $D'D''\pi = - D''D'\pi$, and the above equality entails: 

\vspace{2ex}

$(\mathrm{Id}-\pi)\circ D''D'\pi = -D'\pi \wedge D''\pi = -D'\pi \wedge D''\pi \circ (\mathrm{Id}-\pi),$ 

\vspace{1ex}

\noindent in the sense of currents. This finally gives the following formula:

\vspace{1ex}

$ D''_{\mbox{\tiny det\!\! Q}} D'_{\mbox{\tiny det\!\! Q}}\tau_a = \sum\limits_{j_1, \dots , \, j_q} a_{1 j_1} \dots a_{qj_q} $ \\

\hspace{6ex} $\cdot \sum\limits_k (\mathrm{Id}-\pi)e_{j_1} \wedge \dots  \wedge  (D'\pi\wedge D''\pi)\circ (\mathrm{Id} - \pi) (e_{j_k}) \wedge    \dots \wedge (\mathrm{Id}-\pi)e_{j_q},$

\vspace{1ex}

\noindent and since $\mathrm{Tr}_E(D'\pi \wedge D''\pi) = \mathrm{Tr}_E(\beta \wedge \beta^{\star}),$ we get: \\

$(3. 3)$  \hspace{2ex} $ D''_{\mbox{\tiny det\!\! Q}} D'_{\mbox{\tiny det\!\! Q}}\tau_a =  \mathrm{Tr}_E(\beta \wedge \beta^{\star})\cdot \tau_a .$  \\

\noindent This trivially implies identity $({\star})$. The inequality stated in lemma \ref{Lem:LelPoinc} is thus proved.  

 We will now prove the equality featured by lemma \ref{Lem:LelPoinc}. For $\tau_a$ viewed as an $L^2_1$ section of the bundle $\Lambda^q E$ we get by differentiation: \\

$\begin{array}{lll}

i\partial\bar{\partial}\log (|\tau_a|^2 + \delta^2) & = & i\frac{ \{D_{\Lambda^qE}'\tau_a, \, D_{\Lambda^qE}' \tau_a \}}{|\tau_a|^2 + \delta^2} - i \frac{ \{D'_{\Lambda^qE}\tau_a, \, \tau_a\} \wedge \{\tau_a, \, D'_{\Lambda^qE}\tau_a \}}{(|\tau_a|^2 + \delta^2)^2} \\

\vspace{2ex}

 & - & i\frac{ \{D''_{\Lambda^qE}\tau_a, \, D''_{\Lambda^qE}\tau_a \} }{|\tau_a|^2 + \delta^2} - i \frac{\{\tau_a, \, D''_{\Lambda^qE}\tau_a \}\wedge \{D''_{\Lambda^qE}\tau_a, \, \tau_a \}}{(|\tau_a|^2 + \delta^2)^2} \\

\vspace{2ex}

& + & i\frac{ \{D'_{\Lambda^qE}D''_{\Lambda^qE}\tau_a, \, \tau_a \} }{|\tau_a|^2 + \delta^2} + i \frac{\{\tau_a, \, D''_{\Lambda^qE}D'_{\Lambda^qE}\tau_a \}}{|\tau_a|^2 + \delta^2} \\

\vspace{2ex} 

  & - & i\frac{\{D'_{\Lambda^qE} \tau_a, \, \tau_a\}\wedge \{D''_{\Lambda^qE}\tau_a, \, \tau_a\}}{(|\tau_a|^2 + \delta^2)^2} - i \frac{\{\tau_a, \, D''_{\Lambda^qE}\tau_a\} \wedge \{\tau_a, \, D'_{\Lambda^qE}\tau_a\}}{(|\tau_a|^2 + \delta^2)^2}.

\end{array}$ \\

\noindent The above formulae of $\tau_a$ and $D''_{\Lambda^qE} \tau_a$ imply, for all $k, \, l,$ the identity 

\noindent $\{D''\pi(e_{j_k}), \, (\mathrm{Id}-\pi)e_{j_l} \} = 0$. Indeed, $D''\pi(e_{j_k})$ is an $L^2$ section of $\mathrm{Im}\, \pi,$ and $ \mathrm{Im} \, \pi$ and $\mathrm{Im} \, (\mathrm{Id} - \pi)$ are orthogonal. This readily implies that: \\

$\{D''_{\Lambda^qE}\tau_a, \, \tau_a \} = 0$ \hspace{2ex} and  \hspace{2ex} $\{ \tau_a, \, D''_{\Lambda^qE}\tau_a \} = 0$. \\

\noindent The formula of $i\partial\bar{\partial}\log (|\tau_a|^2 + \delta^2)$ is then reduced to: \\

$\begin{array}{lll}

i\partial\bar{\partial}\log (|\tau_a|^2 + \delta^2) & = & i\frac{ \{D_{\Lambda^qE}'\tau_a, \, D_{\Lambda^qE}' \tau_a \}}{|\tau_a|^2 + \delta^2} - i \frac{ \{D'_{\Lambda^qE}\tau_a, \, \tau_a\} \wedge \{\tau_a, \, D'_{\Lambda^qE}\tau_a \}}{(|\tau_a|^2 + \delta^2)^2} \\

\vspace{2ex}

 & - & i\frac{ \{D''_{\Lambda^qE}\tau_a, \, D''_{\Lambda^qE}\tau_a \} }{|\tau_a|^2 + \delta^2}  \\

\vspace{2ex}

& + & i\frac{ \{D'_{\Lambda^qE}D''_{\Lambda^qE}\tau_a, \, \tau_a \} }{|\tau_a|^2 + \delta^2} + i \frac{\{\tau_a, \, D''_{\Lambda^qE}D'_{\Lambda^qE}\tau_a \}}{|\tau_a|^2 + \delta^2}. \\

\end{array}$ \\

\noindent We will now prove the following three identities: \\

$({\star}{\star})$ \hspace{2ex} $i\{D'_{\Lambda^qE}D''_{\Lambda^qE}\tau_a, \, \tau_a \} = - |\tau_a|^2 \cdot \mathrm{Tr}_E(i\beta \wedge \beta^{\star}), $ \\

$({\star}{\star}{\star})$  \hspace{2ex} $i\{\tau_a, \, D''_{\Lambda^qE}D'_{\Lambda^qE}\tau_a \} =  - |\tau_a|^2  \cdot \mathrm{Tr}_E(i\beta \wedge \beta^{\star}), $ \\
 
$({\star}{\star}{\star}{\star})$ \hspace{2ex} $i\{D''_{\Lambda^qE}\tau_a, \, D''_{\Lambda^qE}\tau_a \} = -  |\tau_a|^2  \cdot \mathrm{Tr}_E(i\beta \wedge \beta^{\star}), $ \\

\noindent and this will prove the equality in the lemma. Let us first prove $({\star}{\star})$. The operator $D'_{\Lambda^qE},$ when applied to the previously obtained formula of $D''_{\Lambda^qE}\tau_a,$ gives: \\

$D'_{\Lambda^qE}D''_{\Lambda^qE}\tau_a = \sum\limits_{j_1, \dots , \, j_q} A_{j_1, \dots , \, j_q},$  \hspace{2ex} where \\

$A_{j_1, \dots , \, j_q}=$ \\

$\begin{array}{lll}

  & = & \sum\limits_{l < k} (\mathrm{Id}-\pi)e_{j_1} \wedge \dots \wedge D'\pi(e_{j_l}) \wedge \dots \wedge D''\pi(e_{j_k}) \wedge \dots \wedge  (\mathrm{Id}-\pi)e_{j_q} \\

\vspace{2ex}

  & - & \sum\limits_k (\mathrm{Id}-\pi)e_{j_1} \wedge \dots \wedge D'D''\pi (e_{j_k}) \wedge \dots \wedge (\mathrm{Id}-\pi)e_{j_q} \\

\vspace{2ex}

  & - & \sum\limits_{k < l} (\mathrm{Id}-\pi)e_{j_1} \wedge \dots \wedge D''\pi(e_{j_k}) \wedge \dots \wedge D'\pi(e_{j_l}) \wedge \dots \wedge  (\mathrm{Id}-\pi)e_{j_q}.

\end{array}$

\noindent Lemma \ref{Lem:distr} (applied to the second sum above) implies that $A_{j_1, \dots , \, j_q}$ is a well-defined current of bidegree $(1, \, 1)$ with values in $E.$

\noindent Since every factor $D''\pi(e_{j_k})$ is orthogonal to every factor $(\mathrm{Id}-\pi)e_{j_l}$ occuring in the expression of $\tau_a,$ we get without difficulty: \\

\noindent $i\{D'_{\Lambda^qE}D''_{\Lambda^qE}\tau_a, \, \tau_a \} = - i\{ \mathrm{Tr}_E (D'\pi \wedge D''\pi) \cdot \tau_a, \, \tau_a \} = - |\tau_a|^2 \cdot \mathrm{Tr}_E (iD'\pi \wedge D''\pi), $ \\

\noindent and this proves $({\star}{\star})$.

\vspace{2ex}

 Let us now prove $({\star}{\star}{\star})$. The computations and arguments are very similar to those of the previous case. These computations show that $D''_{\Lambda^qE}D'_{\Lambda^qE}\tau_a$ is equal, as a $(1, \, 1)$ current, to $D''_{\mbox{\tiny det\!\! Q}}D'_{\mbox{\tiny det\!\! Q}}\tau_a$ plus some terms which have no contribution in the calculation of $i\{\tau_a, \, D''_{\Lambda^qE}D'_{\Lambda^qE}\tau_a \}$. We thus finally get, while taking $(3. 3)$ into account as well: \\

$\begin{array}{lll}

i \{\tau_a, \, D''_{\Lambda^qE}D'_{\Lambda^qE}\tau_a \} & = & i \{\tau_a, \, D''_{\mbox{\tiny det\!\! Q}}D'_{\mbox{\tiny d\'et\!\! Q}}\tau_a \} = i\{\tau_a, \, \mathrm{Tr}_E(\beta \wedge \beta^{\star}) \cdot \tau_a \} \\

  & = & - \{\tau_a, \, \mathrm{Tr}_E(i\beta \wedge \beta^{\star}) \cdot \tau_a \} = - |\tau_a|^2 \cdot \mathrm{Tr}_E(i\beta \wedge \beta^{\star}).

\end{array}$  \\

\noindent This proves $({\star}{\star}{\star})$. It remains to prove $({\star}{\star}{\star}{\star})$. We have already proved that $\{ D''_{\Lambda^qE}\tau_a, \, \tau_a \} =0$ as a $(0, \, 1)$-current with scalar values. If we apply the operator $\partial$ we get: \\

$0 = \partial\{ D''_{\Lambda^qE}\tau_a, \, \tau_a \} = \{D'_{\Lambda^qE}D''_{\Lambda^qE}\tau_a, \, \tau_a \} - \{ D''_{\Lambda^qE}\tau_a, \,  D''_{\Lambda^qE}\tau_a \},$ \\

\noindent which implies, owing to $({\star}{\star})$ : \\

$i\{ D''_{\Lambda^qE}\tau_a, \,  D''_{\Lambda^qE}\tau_a \} = i \{D'_{\Lambda^qE}D''_{\Lambda^qE}\tau_a, \, \tau_a \} = -|\tau_a|^2 \cdot \mathrm{Tr}_E(i\beta\wedge \beta^{\star}).$ \\

\noindent This proves $({\star}{\star}{\star}{\star})$. Relations $({\star})$, $({\star}{\star})$, $({\star}{\star}{\star})$ and $({\star}{\star}{\star}{\star})$ imply the equality stated in lemma \ref{Lem:LelPoinc}. Lemma \ref{Lem:LelPoinc} is thus completely proved.   \hfill $\Box$

\vspace{2ex}

\noindent $\bullet$ {\bf Fourth step} : {\it Construction of the bundle in dimension 1}

 We are now in a position to prove that $\mathrm{Im}\, \pi$ defines almost everywhere a holomorphic subbundle of $E$ in restriction to almost every complex line considered locally in a coordinate patch. As a matter of fact, the proof still works for every complex subspace such that the restriction of the current $\mathrm{Tr}_E(i\beta\wedge \beta^{\star})$ is a well-defined $d$-closed current. Let us fix an arbitrary point $x_0 \in X$ and a trivializing open set $U \ni x_0$ of $E$ contained in a coordinate patch with local coordinates $z=(z_1, \, \dots , \, z_n).$ Let us also fix a complex line $L$ in this coordinate patch such that the restriction of $\mathrm{Tr}_E(i\beta\wedge \beta^{\star})$ to $L$ is a well-defined $(1, \, 1)$-current. This is the case for almost every choice of $L.$ Thanks to corollary \ref{Cor:potentiel}, there exists a subharmonic potential $\varphi = \varphi_L$ on $U\cap L$ such that $i\partial\bar{\partial}\varphi = \mathrm{Tr}_E(i\beta\wedge \beta^{\star})_{|U\cap L}$ (the curvature of $E$ is assumed to be zero according to the first step). Then lemma \ref{Lem:LelPoinc} implies: \\

$\displaystyle i\partial\bar{\partial}\log(|\tau_a|^2 + \delta^2) \geq -\frac{|\tau_a|^2}{|\tau_a|^2 + \delta^2} \, (i\partial\bar{\partial}\varphi) \geq   -i\partial\bar{\partial}\varphi,$  \hspace{2ex} for all $\delta > 0$, \\

\noindent on $U\cap L$, for $i\partial\bar{\partial}\varphi \geq 0$. This shows that the function $\log(|\tau_a|^2 \, e^{\varphi} + \delta^2 \, e^{\varphi})$ is subharmonic on $U\cap L$ for all $\delta > 0,$ and consequently $\log(|\tau_a|^2 \, e^{\varphi})$ is subharmonic on $U\cap L$ as a decreasing limit of subharmonic functions. In particular, the function 

\vspace{1ex}

\hspace{6ex} $\psi=\log(|\tau_a| \, e^{\frac{\varphi}{2}})$ \\

\noindent is subharmonic and not identically $- \infty$ on $U\cap L.$

 Let us now consider a holomorphic function $f: U\cap L \rightarrow \C$ such that $\int_{U\cap L}|f|^2 e^{-2\psi}\, d\lambda < +\infty$, where $d\lambda$ is the Lebesgue measure. The function $\displaystyle |f| \, e^{-\psi} = \frac{|f|}{|\tau_a| \, e^{\frac{\varphi}{2}}}$ is thus $L^2$ on $U\cap L$. In particular, $\displaystyle \frac{f}{\tau_a \, e^{\frac{\varphi}{2}}}$ is an $L^2$ section of $(\det Q)^{-1}$ on $U\cap L$. Since $e^{\frac{\varphi}{2}}$ is subharmonic and, moreover, $L^{\infty}$ on $U\cap L,$ we get that   \\

$\displaystyle \frac{f}{\tau_a} =  e^{\frac{\varphi}{2}} \, \frac{f}{\tau_a \, e^{\frac{\varphi}{2}}}$ \\

\noindent is an $L^2$ section of $(\det Q)^{-1}$ on $U\cap L$. In particular, it defines a distribution and the expression $D''\bigg (\frac{f}{\tau_a} \bigg )$ is well-defined in the sense of distributions. We have thus obtained the regularity needed for the application of the operator $D''$ (as explained at the beginning of the third step).

 Now the arguments enabling us to conclude are purely formal. Indeed, $D''\bigg (\frac{f}{\tau_a} \bigg) = 0$ at all points where this is well-defined. The bundle morphism $v$ defined by $(3. 2)$ can then be redefined on $U\cap L$ as: \\

$\displaystyle \Lambda^{q+1}E  \stackrel{v}{\longrightarrow} E$, \hspace{2ex} $e_I \mapsto \frac{f\, u(e_I)}{\tau_a},$ \\

\noindent for all multiindex $I$ such that $|I|=q+1$. Since $u(e_I)$ is a $D''$-closed $L^2$ section of $E \otimes \det Q$ on $U\cap L,$ we get that $\frac{f\, u(e_I)}{\tau_a} \in L^1(U\cap L, \, E)$ and \\

 $\displaystyle D''\bigg (\frac{f\, u(e_I)}{\tau_a} \bigg ) = D''\bigg ( \frac{f}{\tau_a}\bigg) \, u(e_I) +  \frac{f}{\tau_a} \, D''u(e_I) = 0,$ \\

\noindent for all $I$. Consequently, the $L^2$ bundle defined by $F=\mathrm{Im}\, v = \mathrm{Im}\, \pi$ is locally generated by its local meromorphic sections $\displaystyle \frac{f\, u(e_I)}{\tau_a}$ on almost every complex line $L$ contained in a coordinate patch.

\vspace{2ex}

\vspace{2ex}

\noindent $\bullet$ {\bf Fifth step} : {\it application of a theorem of Shiffman's}

 This step would be superfluous if we were able to prove that the current $\mathrm{Tr}_E(i\beta\wedge \beta^{\star})$ is $d$-closed (see the explanation at the beginning of the previous step) . Recall that $p$ is the rank almost everywhere of $\pi.$ Let us now consider the following map relative to the trivializing open set $U$ of $E$: \\

$U \ni x \stackrel{\Phi}{\longmapsto} G(p, \, r)$  \\

\noindent definied almost everywhere as $\Phi(x) = \mathrm{Im}\, \pi_x.$ This is a $p$-dimensional vector subspace of $E_x$, and it can therefore be viewed as an element in the Grassmannian $G(p, \, r)$ of $p$-dimensional vector subspaces of $\C^r$. As the Grassmannian is a projective manifold, there exists an isometric embedding of $G(p, \, r)$ into the complex projective space $\mathbb{P}^K,$ which can in its turn be embedded into a Euclidian space $\R^N$. The vector-valued map \\

$\Phi = (\Phi_1, \, \dots , \, \Phi_N) : U \rightarrow G(p, \, r) \hookrightarrow \R^N$ \\

\noindent is $L^2_1$ for it is defined by  $\pi$ which is assumed to be $L^2_1.$ What we have proved above amounts to the component $\Phi_j : U\cap L \rightarrow \R$ of $\Phi$ being meromorphic almost everywhere for almost all complex line $L$ and all $j=1, \, \dots , \, N$. The following Hartogs-type theorem is due to B.Shiffman (see [Shi86], corollary 2, page 240). It states that a measurable function which is separately meromorphic almost everywhere is in fact meromorphic almost everywhere. \\

\noindent {\bf Theorem (Shiffman, 1986).} {\it Let $\Delta$ be the unit disc of $\C$ and let $f : \Delta^n \longrightarrow \C$ be a measurable function such that for all $1 \leq j \leq n$ and almost all $(z_1,  \dots , \, \hat{z}_j, \, \dots , \, z_n) \in \Delta^{n-1}$, the map $\Delta \ni z_j \mapsto f(z_1, \dots , \, z_n)$ is equal almost everywhere to a meromorphic function on $\Delta$. Then $f$ is equal almost everywhere to a meromorphic function.} \\

  It is noteworthy that the hypotheses of this theorem of Shiffman's are quite loose. The function $f$ is merely assumed to be measurable and meromorphic almost everywhere along almost all directions parallel to the coordinate axes. Our above-defined functions $\Phi_j$ satisfy much stronger hypotheses. They are not only measurable but also $L^2_1$. They are meromorphic almost everywhere along almost all directions as well.

 This result implies that the components $\Phi_j$ of $\Phi$ are meromorphic almost everywhere. The map $\Phi$ is thus meromorphic almost everywhere. Since every meromorphic map is holomorphic outside an analytic subset of codimension $\geq 2$, we get that $F= \mathrm{Im}\, \pi$ is a holomorphic subbundle of $E$ outside an analytic subset $S\subset X$ of codimension $\geq 2$.

\vspace{2ex}

\noindent {\bf Acknowledgements.} {\it I am grateful to my thesis supervisor Jean-Pierre Demailly for his unflinching support and great scientific expertise.}

\vspace{3ex}

\begin{center}

{\bf References}

\end{center}

\vspace{2ex}

\noindent [Dem 97] \, J.-P. Demailly --- {\it Complex Analytic and Algebraic Geometry}---http://www-fourier.ujf-grenoble.fr/~demailly/manuscripts/agbook.ps.gz

\vspace{1ex}

\noindent [Gri69] \, P. A. Griffiths --- {\it Hermitian Differential Geometry, Chern Classes and Positive Vector Bundles}--- Global Analysis, papers in honour of K. Kodaira, Princeton Univ. Press, Princeton, 1969, 181-251.

\vspace{1ex}

\noindent [Kob87] \, S. Kobayashi --- {\it Differential Geometry of Complex Vector Bundles} --- Princeton University Press, 1987.

\vspace{1ex}

\noindent [LT95] \, M. L\"ubke, A. Teleman --- {\it The Kobayashi-Hitchin Correspondence}--- World Scientific, 1995.

\vspace{1ex}

\noindent [Pop03] \, D.Popovici --- {\it Quelques applications des m\'ethodes effectives en g\'eom\'etrie analytique} --- PhD thesis, Universit\'e de Grenoble, 2003.

\vspace{1ex}

\noindent [Shi86] \, B. Shiffman --- {\it Complete Characterization of Holomorphic Chains of Codimension One}--- Math. Ann., {\bf 274}, (1986) P. 233-256.

\vspace{1ex}

\noindent [Sib85] \, N. Sibony --- {\it Quelques probl\`emes de prolongement de courants en analyse complexe} --- Duke Math. J. , {\bf 52 }, No. {\bf 1} (1985), P. 157 - 197.

\vspace{1ex}

\noindent [UY 86]\, K. Uhlenbeck, S.T. Yau --- {\it On the Existence of Hermitian-Yang-Mills Connections in Stable Vector Bundles} --- Communications in Pure and Applied Mathematics, Vol. {\bf XXXIX}, 1986, Supplement, pp. S257 - S293.

\vspace{1ex}

\noindent [UY 89] \, K. Uhlenbeck, S.T. Yau ---{\it A Note on Our Previous Paper: On the Existence of Hermitian-Yang-Mills Connections in Stable Vector Bundles} --- Communications on Pure and Applied Mathematics, Vol. {\bf XLII}, pp. 703 - 707.

\vspace{3ex}

\noindent Dan Popovici

\noindent Universit\'e de Grenoble I, Institut Fourier, BP74, 38402 Saint-Martin d'H\`eres Cedex, France

\noindent E-mail:popovici@ujf-grenoble.fr

\vspace{3ex}

\noindent Universit\'e de Paris-Sud, Laboratoire de Math\'ematiques, B\^at. 425, 91 405 Orsay Cedex, France

\noindent E-mail:Dan.Popovici@math.u-psud.fr

\end{document}